\documentclass[12]{amsart}
\usepackage{graphicx}
\usepackage{amsmath}
\usepackage{authblk}
\usepackage{amsthm}
\usepackage{amssymb}
\usepackage{xcolor}

\usepackage{tikz}
\usepackage{pgfplots}
\pgfplotsset{compat=1.13}
\usetikzlibrary{intersections, patterns, pgfplots.fillbetween}

\usepackage{braids}
\usepackage{amsaddr}

\usepackage{geometry}
\newtheorem{theorem}{Theorem}[section]
\newtheorem{corollary}{Corollary}[theorem]
\newtheorem{lemma}[theorem]{Lemma}
\newtheorem{example}[theorem]{Example}

\theoremstyle{definition}
\newtheorem{definition}{Definition}[section]

\theoremstyle{remark}

\begin{document}

\title{A Canonical Form for Max Plus Symmetric Matrices and Applications}

\author{ Himadri Mukherjee$^1$, Askar Ali M$^2$}
\address{$^{1,2}$ Department of Mathematics,\\ BITS Pilani K. K. Birla Goa Campus, Goa, India}
\email{$^1$ himadrim@goa.bits-pilani.ac.in, $^2$ p20190037@goa.bits-pilani.ac.in}
\footnote{AMS Classification 2020. Primary: 15A24, 16T25}

\maketitle

\noindent

\begin{abstract}
We develop a canonical form for congruence of max plus symmetric matrices. We use the same canonical form to get results in the generalized eigenvector problem. We have also utilized the canonical form to find all symmetric matrices that commute with a given symmetric matrix.
\end{abstract}



\section{Introduction} 
Let $\mathbb{T}$ denote the tropical semi-field, defined as ($\mathbb{T} = \mathbb{R} \cup \{-\infty\}, \oplus, \otimes$), such that, for $a, b \in \mathbb{T},\  a \oplus b := max \{a,b\}$ and $a \otimes b := a+b$, the usual sum in real line. As a convention, for the additional element $-\infty$, which will be denoted by $\varepsilon$ more often in this paper, we have, for any $a \in \mathbb{R},\  a \oplus \varepsilon = a$, and $a \otimes \varepsilon = \varepsilon$. \\
In this article, we study the conjugate action of tropical matrices. i.e., a matrix $P$ acts on $A$ as $P^tAP$. Recall that two matrices $A, B \in M_n(\mathbb{T})$ are said to be congruent if there exists a $P \in Gl_n(\mathbb{T})$, such that $B = P^tAP$.  In that case, we denote $A \sim B$. Note that the only class of invertible matrices in max-plus algebra are generalised permutation matrices \cite{economics-note}. This congruence relation gives an equivalence relation in matrix algebra. Hence, it is natural to think of a canonical form for such equivalent classes. In classical linear algebra, numerous studies have been conducted on these conjugate actions and resulting monoids and orbits (see \cite{teran,conjugate_1,conjugate_2,conjugate_3,teran_2}). In \cite{canonical_general}, H. W. Turnbull et al., proved that any complex square matrix is congruent to a direct sum of blocks of at most six canonical forms. The most interesting canonical forms came from the class of symmetric and skew-symmetric matrices. Real symmetric matrices are completely classified by Sylvester's law of inertia \cite{sylvester_inertia}. For the class of skew-symmetric matrices, we have the following canonical structure \cite{darboux1}. 
\[\begin{pmatrix}
    0&\mu_1&&&&&&\\
    -\mu_1&0&&&&&&\\
    &&0&\mu_2&&&&\\
    &&-\mu_2&0&&&&\\
    &&&&\ddots&&&\\
    &&&&&0&&\\
     &&&&&&\ddots&\\
     &&&&&&&0
    
\end{pmatrix}\]
Such studies of congruent matrices come with plenty of applications towards the fields of eigenvalue-eigenvector problems, matrix pencils and bilinear forms, because of their preserving nature and application in numerical methods \cite{congruence_2,conjugate_3,pencils}.
\\
The story is not different in tropical algebra as well. Recently, in \cite{crypto}, authors have found the application of tropical congruent transformation of symmetric matrices in cryptography. In this article, we provide a canonical form for symmetric matrices through the conjugate action in a pseudo-diagonal matrix form. For certain symmetric matrices $A$ (pseudo-diagonalizable), through a  $2 \times 2$ minor, which is defined by $\begin{vmatrix}
    a_{ij}&a_{ik}\\
    a_{lj}&a_{lk}
\end{vmatrix}_{T} = (a_{ij}\otimes a_{lk}) \otimes (a_{lj} \otimes a_{ik})^{-1}$, we have, for some generalised permutation matrix $P$ such that,  \[P^tAP = \begin{pmatrix}
        d_1 & 0 &0 &\dots &0\\
        0&d_2&0&\dots & 0\\
       \vdots\\
        0&0&0&\dots & d_n
    \end{pmatrix}, \text{ where } d_i+d_j = \begin{vmatrix}
         a_{ii}&a_{ij}\\
    a_{ij}&a_{jj}
    \end{vmatrix}_{T}.  \]  
With this canonical form, we have found various applications in solving generalised eigenvalue problems to finding all commuting symmetric pairs. The problem of solving generalised eigenvalue \cite{Binding.gen.eigen}, i.e., for a non-trivial vector $x$, and a scalar $\lambda$, $A \otimes x = \lambda \otimes B \otimes x$, draws a parallel to the theory of matrix pencils, in the max-plus algebra. The generalised eigenvalue problem has found multiple applications in fields such as optimization and synchronization in multi-machine interactive production processes and various other practical scenarios (see \cite{thebook} for details).\\
  Defining skew-symmetric matrices in a classical way is not possible in max-plus algebra, as the `addition' operation does not have an inverse. However, defining skew-symmetric matrices in a way that enables us to split any matrix into the sum of a symmetric matrix and a skew-symmetric matrix brings a lot of practical usages. We have defined skew-symmetric matrices in the following way. $A = [a_{ij}]$ is skew-symmetric if it holds the following. 
  \[a_{ij} = 
    \varepsilon, \text{ whenever } a_{ji} \neq \varepsilon   \]
We have given a few canonical structures for such skew-symmetric matrices. However, classifying all of them is difficult, as it is equivalent to classifying all simple directed graphs. 
  
\subsection{Prerequisites}

\begin{itemize}
    \item Both $N$ and $[n]$ represents the set $\{1,2, \cdots,n\}$, for appropriate $n$.
    \item  For any $a, b \in \mathbb{T}$, $a \oplus' b := min \{a,b\}$, and for matrices $A,B$, the corresponding product $(A \otimes' B)_{ij} = \underset{k}{min} \{a_{ik}+b_{kj}\}$.
    \item $\mathbb{T^*}$ represents the set $\mathbb{T} - \{\varepsilon \}$.
    \item For two matrices $A,B$, $A \boxtimes B$ represents the Kronecker tensor product of $A$ and $B$ in max-plus algebra.
\end{itemize}
\begin{definition}
    Two matrices $A, B$ are called congruent if there exists a $P \in Gl_n(\mathbb{T}) $, such that $A = P^t B P$. In that case, we denote $A \sim B$. 
\end{definition}
Many of the theories discussed in this paper are closely related to one-sided max-linear systems or two-sided max-linear systems of equations. We make use of the existing theory in this regard. One such key tool is referred to here. Readers may go through \cite{thebook,one-sided-cunninghame} for detailed understanding.\\
Consider the one-sided max-linear system, 
\begin{equation}\label{one-sided-system}
    Ax=b
\end{equation}
where, $A \in M_n(\mathbb{T})$ and $b \in \mathbb{T}^n$. Then, we define $\Bar{x} = (\Bar{x_1},\Bar{x_2}, \cdots, \Bar{x_n})^t$, where $\Bar{x_j} = (\underset{i \in N}{max}~a_{ij}\otimes b_i^{-1})^{-1}$. In other words, 
\[ \Bar{x} = A^{\sharp} \otimes' b,~ \text{ where, } A^ {\sharp} = [-a_{ij}].\]
Now, define the set $N_j(A,b)$ (or simply $N_j$), as follows.
\[N_j (A,b) = \{i \in N| ~ \Bar{x_j} = b_i \otimes a_{ij}^{-1}\}\]
\begin{theorem}[Cunninghame-Green RA,\cite{one-sided-cunninghame}]
  Let $A \in M_n (\mathbb{T})$, and $b \in \mathbb{T^*}^n$. Then the following statements are equivalent:
  \begin{itemize}
      \item[(a)] The system in \eqref{one-sided-system} is solvable,
      \item[(b)] $\Bar{x}$ is a solution to \eqref{one-sided-system},
      \item[(c)] $\cup _{i \in N}N_j = N$. 
  \end{itemize}  
\end{theorem}
\section{Conjugate action and associated monoids}
Let $M= M_n(\mathbb{T})$ be the set of all $n \times n$ tropical matrices. Then, note that $M$ is a monoid. i.e., $X, Y \in M \implies XY \in M$ and $I \in M$. \\
A map $\phi : M \times M \rightarrow M$ is called a monoid action if it satisfies the following:
\begin{itemize}
    \item[i)] $I*m = m,~ \forall m \in M$
    \item[ii)] $X*(Y*m) = (XY)*m$, 
\end{itemize}
where, $\phi(X,m) = X*m$.\\
\begin{example}
  \begin{itemize}
      \item[i)] $X*m = Xm$
      \item[ii)] $X*m = XmX^t$
      \item[iii)] $X*m = mX^t$
  \end{itemize}  
\end{example}
Given $m \in M$, $G_m := \{X | X*m = m\}$, are sub-monoids of $M$.  For a given matrix $A \in M$, $Sol_A := \{X | X^t A X = A\}$. The study of such monoids finds interesting applications. Its classical linear algebra analogue can be found in \cite{HMGS}.
\begin{lemma}
    Let $A = \begin{pmatrix}
        \varepsilon&0&\varepsilon\\
        \varepsilon&\varepsilon&0\\
        \varepsilon&\varepsilon&\varepsilon
    \end{pmatrix}$, then $Sol_A = \{ \begin{pmatrix}
       a&\varepsilon&\varepsilon\\
        \varepsilon&-a&\varepsilon\\
        \varepsilon&\varepsilon&a 
    \end{pmatrix} |~ a \in \mathbb{T} \}$.
\end{lemma}
\begin{proof}
    Let $X = \begin{pmatrix}
        a&b&c\\
        d&e&f\\
        g&h&k
    \end{pmatrix}$. Then,
    \[ X^t AX = \begin{pmatrix}
        ad\oplus dg& ae\oplus dh& af\oplus dk\\
        bd \oplus eg & be \oplus eh & bf \oplus ek\\
        cd \oplus fg & ce \oplus fh & cf \oplus fk
    \end{pmatrix}\]
    Now, equating with $A$ gives the above solution. 
\end{proof}
\begin{lemma}
    Let $A = \begin{pmatrix}
        \varepsilon&\varepsilon& \dots &0\\
        \varepsilon&\varepsilon& \dots &\varepsilon\\
       \vdots\\
        \varepsilon&\varepsilon& \dots &\varepsilon
    \end{pmatrix}$. Then, $Sol_A = \{ \begin{pmatrix}
      a&\varepsilon& \dots &\varepsilon \\
        &D&  & \\
        \varepsilon& \dots &\varepsilon&-a
    \end{pmatrix} | D \in M_{(n-2)\times(n)}(\mathbb{T}), a \in \mathbb{T^*}\}$.
\end{lemma}
\begin{proof}
  Let $X = \begin{pmatrix}
      x_{11}&x_{12}&\cdots&x_{1n}\\
      x_{21}&x_{22}&\cdots&x_{2n}\\
      \vdots\\
      x_{n1}&x_{n2}&\cdots&x_{nn}\\
  \end{pmatrix}$. Then, upon equating $X^tAX = A$, we get,
  \[x_{11} = -x_{nn} \neq \varepsilon,\]
  \[\begin{cases}
      x_{ni} = \varepsilon, \text{ for } 1\leq i <n\\
      x_{1i} = \varepsilon, \text{ for } 1< i \leq n
  \end{cases}\]
  For the rest of the variables $x_{ij}$, irrespective of their values, the equation is satisfied. Hence, they can take any arbitrary values. 
\end{proof}
\begin{lemma}
\begin{itemize}
    \item[(a)]  $Sol_I = \Theta (n)$, where $\Theta(n)$ is collection of all $n \times n$ orthogonal matrices.
    \item[(b)]  Let $A = diag (a_1,a_2, \dots , a_n)$, such that $\varepsilon < a_i$, for $i = 1,2, \dots , n$. Then, $Sol_A = ~ \{X \in Gl_n(\mathbb{T})| \text{ whenever } (X)_{ij} \neq \varepsilon, (X)_{ij} = \frac{a_j- a_i}{2}\}.$
\end{itemize}
   
\end{lemma}
\begin{proof}
\begin{itemize}
    \item[(a)] Let $A = I$, then $X^t A X = A$, implies $X^tX = I$.\\ 
    This gives, $Sol_I = \Theta$.
    \item[(b)] Since $A$ is invertible, $X \in Sol_A$ has to be invertible. Hence, we can assume $X$ to be a generalised permutation matrix. Then, the Rest of the proof is straightforward from direct calculations. 
    \end{itemize}
\end{proof}
\begin{example}
     Let $A$ be the diagonal matrix $\begin{pmatrix}
        a_1 & \varepsilon\\
        \varepsilon & a_2
    \end{pmatrix}$ , then $Sol_A = \{ I, \begin{pmatrix}
     \varepsilon &  \frac{a_2 - a_1}{2}\\
     \frac{a_1 - a_2}{2} & \varepsilon
    \end{pmatrix} \} \simeq Z_2 $.
\end{example}
\begin{lemma}
   If $A \sim B$, then $Sol_A \sim Sol_B$. 
\end{lemma}
\begin{proof}
    Let $A \sim B$. Then, there exists a $P \in Gl_n(\mathbb{T})$, such that $A = P^tBP$. Now, if $X \in Sol_A$, then, 
    \begin{equation*}
        \begin{aligned}
            X^tAX &= A\\
            \implies P^tX^tPP^tAPP^tXP &= P^tAP\\
            \implies Y^tBY &= B,~ \text{ where } Y=P^tXP.
        \end{aligned}
    \end{equation*}
This gives, $Sol_A \sim Sol_B$.
\end{proof}
With the semi-ring structure of tropical algebra, there are various notions of ranks and regularity that do not coincide. For example, Tropical rank, Gondran-Minoux rank, column rank, row rank, Kapranov rank, Barvinoc rank and many more (\cite{rank,akian-rank, ultimate-rank,butkovic-rank}). With respect to these rank notions, we can compare the rank of $A$ and the rank of $X \in Sol_A$. 
\begin{definition}
    A matrix $A = [a_{ij}] \in M_n (\mathbb{T})$, is called non-singular, if there exists a unique permutation $\sigma \in S_n$, such that its permanent, $maper(A) = a_{1\sigma(1)}+a_{2\sigma(2)}+ \cdots +a_{n\sigma(n)} $. 
\end{definition}
\begin{definition}
    The tropical rank $(rk_{tr})$ of a matrix $A \in M_n(\mathbb{T})$ is defined as the maximal $r$, for which there exists an $r \times r$ non-singular sub-matrix of $A$. If $A$ has tropical rank $n$, then we say $A$ is strongly regular.
\end{definition}
\begin{definition}
    Let $A \in M_n(\mathbb{T})$, and $A_n \subset S_n$ denotes the set of even permutations. Then, 
    \[|A|^+ := \underset{\sigma \in A_n}{Max}~ a_{1\sigma(1)}+a_{2\sigma(2)}+ \cdots +a_{n\sigma(n)}\]
     \[|A|^- := \underset{\sigma \in S_n - A_n}{Max}~ a_{1\sigma(1)}+a_{2\sigma(2)}+ \cdots +a_{n\sigma(n)}\]
\end{definition}
\begin{definition}
    The determinantal rank of a matrix $A \in M_n(\mathbb{T})$, denoted by $rk_{det}$, is the maximal $r$, for which there exists an $r \times r$ sub-matrix $A'$ of $A$, such that $|A|^+ \neq |A|^-$. We say, $A$ is determinantaly regular, if $rk_{det}(A)= n$. 
\end{definition}

\begin{theorem}
    Let $A \in M_n(\mathbb{T})$, and $X \in Sol_A$. Then, we have the following.
    \begin{itemize}
        \item [(a)] If $A$ is strongly regular, so is $X$.  
        \item[(b)] If $A$ is determinantaly regular, so is $X$. 
    \end{itemize}
\end{theorem}
\begin{proof}
    The proof follows from the fact that (see \cite{akian-rank} for details),
    \begin{align*}
         rk_{tr}(AB) &\leq ~ min~\{rk_{tr}(A), rk_{tr}(B)\}, \\ 
        rk_{det}(AB) &\leq~ min~\{rk_{det}(A), rk_{det}(B)\}.
    \end{align*}  
\end{proof}
\section{Canonical Forms of Symmetric and Skew Symmetric matrices}
The question of finding a canonical form for the congruence of tropical matrices can be thought of as a problem of finding the quotient of a certain vector space with respect to an action of a group. Let us denote the set of $n \times n$ tropical symmetric matrices by $\Delta_n(\mathbb{T})$. Let us take the action of the general linear group $GL_n(\mathbb{T})$ on $\Delta_n(\mathbb{T})$ by the congruence action, namely an invertible matrix $P$ acts on a symmetric matrix $X$ by sending $X$ to $P^tXP$. The question of finding a canonical form is the same as finding a set of coset representatives for $\Delta_n(\mathbb{T})/ GL_n(\mathbb{T})$. 
This section describes a canonical structure of tropical symmetric and skew-symmetric matrices. As the $``\oplus"$ operation does not have an inverse, it is not possible to define skew-symmetric matrices in a classical way. However, the following definition gives many practical uses in max-plus algebra. 
\begin{definition}\label{skew-sym}
    A matrix $A \in M_n(\mathbb{T})$ is said to be skew-symmetric, if $a_{ij} = 
    \varepsilon, \text{ whenever } a_{ji} \neq \varepsilon  $
\end{definition}
With the above definition, we can write any tropical square matrix as a sum of symmetric and skew-symmetric matrices. i.e., for any $A \in M_n(\mathbb{T})$, we have a symmetric matrix $S$ and skew symmetric matrix $U$, such that $A = S \oplus U$. This can be achieved in the following way.\\
Set the diagonal entries as, $u_{ii} = \varepsilon$ and $s_{ii} = a_{ii}$. For the remaining entries,
\begin{equation*}
   \text{ if } a_{ij} > a_{ji}, \text{ set } \begin{cases}
    u_{ij} = a_{ij},~ u_{ji}= \varepsilon\\
    s_{ij} = s_{ji} = a_{ji}
    \end{cases}, 
    \text{ and if } a_{ij} = a_{ji}, \text{ set } \begin{cases}
        s_{ij} = s_{ji} = a_{ij}\\
        u_{ij} = u_{ji} = \varepsilon
    \end{cases}  
\end{equation*}
    If $a_{ji} > a_{ij}$, then swap the role of $i,j$ in the above setting. Note that if there are no symmetric entries (i.e., $a_{ij} \neq a_{ji}$, for any $i,j$), then the splitting of the matrix $A$ into symmetric and skew-symmetric is unique.
 \begin{example}
   Let $A = \begin{pmatrix}
       4&6\\
       7&9
   \end{pmatrix}$ . Take, $S = \begin{pmatrix}
       4&6\\
       6&9
   \end{pmatrix}$  and $U = \begin{pmatrix}
       \varepsilon&\varepsilon\\
       7&\varepsilon
   \end{pmatrix}$. Then, we have $A = S \oplus U$, and there are no other symmetric and skew-symmetric matrices that sum up to $A$.  
 \end{example}
 \textbf{Remark:} Let $A$ be a symmetric matrix (or skew-symmetric). Then, for any generalised permutation matrix $P$, we have $P^tAP$ is also symmetric (or skew-symmetric).\\
 To derive the canonical forms of finite symmetric matrices, we use the following pseudo-determinant on $2 \times 2$ minors. For a matrix $A = [a_{ij}]$,  
 \begin{equation}\label{pseudo-det}
     \begin{vmatrix}
 a_{ik}&a_{il}\\
 a_{jk}&a_{jl}
 \end{vmatrix}_T = \Big(\begin{vmatrix}
 a_{ik}&a_{il}\\
 a_{jk}&a_{jl}
 \end{vmatrix}^+\Big)\otimes \Big(\begin{vmatrix}
  a_{ik}&a_{il}\\
 a_{jk}&a_{jl}
 \end{vmatrix}^-\Big)^{-1} =  (a_{ik}\otimes a_{jl}) \otimes (a_{il}\otimes a_{jk})^{-1} 
 \end{equation}

Let us denote the set of multi-indices $\{ (i_1,i_2, \ldots, i_d) \, \mid \, 1 \leq i_1 < i_2 \cdots < i_d \leq n\}$ by $I_{d,n}$. With this notation in place, for two multi-indices $I,J \in I_{d,n}$, $(I\mid J)$ stands for the ``co-factor" of a matrix, which is the pseudo determinant of the sub-matrix with rows in $I$ and columns in $J$. In other words, consider $I_{2,n} = \{\underline{i}= (i_1,i_2): 1 \leq i_1 < i_2 \leq n \}$. Explicitly, $(\mbox{ rows } \mid \mbox{ columns }) = (i,j\mid l,k)$ from $I_{2,n} \times I_{2,n}$ cofactor of $A$, is denoted by $(i,j\mid  l,k) \rhd A$, and is equal to   $(a_{ik}\otimes a_{jl}) \otimes (a_{il}\otimes a_{jk})^{-1}$. In terms of usual algebra, the same quantity can be expressed as $ (a_{ik} + a_{jl}) - (a_{il}+ a_{jk})$. This pseudo-determinant has found usage in many recent works, notably in \cite{matrix_root}. With respect to the lexicographic ordering on $I_{2,n}$, we can define a partial ordering on $I_{2,n} \times I_{2,n}$, as $(\underline{i},\underline{j}) \leq (\underline{k},\underline{l})$, iff $\underline{i} \leq \underline{k}$, and if $\underline{i} = \underline{k}$, then $\underline{j} \leq \underline{l}$. 
 One can observe that this pseudo determinant is invariant under scalar multiplication (or scaling). i.e.,  $( \underline{i}\mid\underline{j}) \rhd A =  (\underline{i}\mid\underline{j}) \rhd \lambda A $, for any $\lambda \in \mathbb{T}^*$.  
 \begin{definition}
     Let $\Psi$ be the map from $\mathrm{M}_n(\mathbb{T}^*)$ to $ \mathrm{M}_{\binom{n}{2}}(\mathbb{T}^*)$, which maps $ A \mapsto \Bar{A}$,
     where $\Bar{A}$ is obtained by applying the multi-index operator $(i,j\mid l,k)$ on $A$ in the increasing order. For convenience, we use this multi-index for the elements of $ \mathrm{M}_{\binom{n}{2}}(\mathbb{T}^*)$ eg. the $((i,j),(l,k))-th$ entry of $\Bar{A}$ is  $(i,j\mid  l,k) \rhd A$. An $S_n$ action on $\Bar{A}$ is defined as $((i,j),(l,k))-th$ entry of $\Bar{A}$ goes to $((\sigma(i), \sigma(j)),(\sigma(l),\sigma(k)))-th$ entry, for some $\sigma \in S_n$.
 \end{definition}

 The map $\Psi$ takes the line $\lambda A$ to the point $\Bar{A}$. Also, note that multiple lines can be mapped to the same $\Bar{A}$.\\
 Interestingly, this multi-index operator acts similarly to $``Pl\ddot{u}cker$ relation". For $i_1 < i_2 < i_3$, $1 \leq j,k \leq n$, it follows,
 \begin{equation}\label{Plucker-relation.1}  
 (i_1,i_3\mid j,k)  -   (i_1,i_2\mid j,k) = (i_2,i_3\mid j,k) 
 \end{equation}
 \begin{equation}\label{Plucker-relation.2}
    (j,k \mid i_1,i_3)  -   (j,k \mid i_1,i_2) = (j,k\mid i_1,i_2) 
 \end{equation}
 \begin{example}
     Let $A= \begin{pmatrix}
         1&3&4&5\\
         3&2&6&7\\
         4&6&-3&8\\
         5&7&8&6
     \end{pmatrix}$. Then $\Bar{A} = \Psi(A) = \begin{pmatrix}
         -3&0&0&3&3&0\\
         0&-10&0&-10&0&10\\
         0&0&-3&0&-3&-3\\
         3&-10&0&-13&-3&10\\
         3&0&-3&-3&-6&-3\\
         0&10&-3&10&-3&-13
     \end{pmatrix}$.\\
Let $R_i, r_i$ represent the $i^{th}$ row of $\Bar{A}$ and  $A$ respectively. Then, $R_1$ is formed using $r_1$ and $r_2$; $R_2$ using $r_1$ and $r_3$; and $R_3$ using $r_1$ and $r_4$, and so on. One can observe that $R_2 - R_1 = R_4$. This corresponds to the $Pl\ddot{u}cker$ like relation,
 \[ (1,3\mid j,k) -  (1,2\mid j,k) = (2,3\mid j,k). \]
 Similarly, we have $R_3 - R_1 = R_5$, $R_3 - R_2 = R_6$ and so on. We can observe similar relations in the columns as well.
 \end{example}
  Mapping an $n \times n$ matrix to an $\binom{n}{2} \times \binom{n}{2}$ matrix seems to increase the complexity of any arithmetic with them. But in reality, any $\Bar{A}$ has an $(n-1) \times (n-1)$ principal block in it, which generates the remaining blocks of the matrix. Any matrix that follows $Pl\ddot{u}cker$ like relations as in \eqref{Plucker-relation.1}-\eqref{Plucker-relation.2}, the left-top $(n-1) \times (n-1)$ corner block act as principal block. Let us demonstrate this using the above example. We have, $\Bar{A} = \begin{pmatrix}
      \Bar{A}_1&\Bar{A}_2\\
      \Bar{A}_3&\Bar{A}_4
  \end{pmatrix} = \left(\begin{array}{c c c| c c c} 
	-3&0&0&3&3&0\\
         0&-10&0&-10&0&10\\
         0&0&-3&0&-3&-3\\
         \hline
         3&-10&0&-13&-3&10\\
         3&0&-3&-3&-6&-3\\
         0&10&-3&10&-3&-13
	
\end{array}\right) 
$, the left-top corner $3 \times 3$ block $\Bar{A}_1$ generate remaining. First, generate the right-top corner $3 \times 3$ block $\Bar{A}_2$ by taking the differences of columns in  $\Bar{A}_1$ (with respect to \eqref{Plucker-relation.2}). Now, to generate $\Bar{A}_3$ and $\Bar{A}_4$, take the differences of rows of the combined matrix $(\Bar{A}_1~\Bar{A}_2)$ (with respect to \eqref{Plucker-relation.1}).    
 \begin{theorem}\label{lemma_action}
     Let $A \sim B$, such that for some generalized permutation matrix $P$,  $P^tAP= B$, and $\sigma \in S_n$ be the corresponding permutation to $P$. Then,
     \[  (\underline{i}\mid \underline{j}) \rhd A = ( \sigma(\underline{i})\mid  \sigma(\underline{j})) 
          \rhd  B \]
    Where, $\sigma (\underline{i}) = (\sigma(i_1), \sigma(i_2))$. The converse is true if $A$ and $B$ are finite symmetric matrices.
 \end{theorem}
 \begin{proof}
     Let $P^tAP = B$, and $\sigma$ be the corresponding permutation associated with $P$. Under the action of $P$, an entry $a_{i_{m} i_{k}}$ in $A$, goes to $b_{\sigma(i_{m})\sigma(i_{k})}$ in $B$. With this action, each index $i_m$ carries a weight $p_{i_m}$, where $p_{i_m}$ is a finite entry in $P$. i.e., under the action of $P$, $a_{i_{1} i_{2}}$ in $A$ goes to $b_{\sigma(i_{1})\sigma(i_{2})} = a_{i_{1} i_{2}} \otimes p_{i_1} \otimes p_{i_2}$. Hence,
     \begin{align*}
          (\sigma(\underline{i})\mid  \sigma(\underline{j}) ) 
          \rhd  B  &= (b_{\sigma(i_1)\sigma(j_1)}+ b_{\sigma(i_2)\sigma(j_2)}) - (b_{\sigma(i_1)\sigma(j_2)} + b_{\sigma(i_2)\sigma(j_1)})\\
          &= \Big((a_{i_1j_1} + p_{i_1} + p_{j_1})+ (a_{i_2j_2}\otimes p_{i_2} \otimes p_{j_2})\Big) -  \Big((a_{i_1j_2} + p_{i_1} + p_{j_2})+ (a_{i_2j_1}\otimes p_{i_2} \otimes p_{j_1})\Big)\\
          &= (a_{i_1j_1} + a_{i_2j_2}) - (a_{i_1j_2}+ a_{i_2j_1})\\
          &=  (\underline{i}\mid \underline{j}) \rhd A
     \end{align*} 
Now, conversely assume that, $A$ and $B$ are finite symmetric matrices, and for some $\sigma \in S_n$, $ (\underline{i}\mid \underline{j}) \rhd A = (\sigma(\underline{i})\mid  \sigma(\underline{j}))\rhd  B$, for every $\underline{i},\underline{j} \in I_{2,n}$. Let $P_0$ be the permutation matrix corresponding to $\sigma$. Define $A' := P_0^tAP_0$. Then we have, $a'_{\sigma(i_m)\sigma(i_k)} = a_{i_mi_k}$. For convenience, take $i_m' = \sigma (i_m)$. Then, define the matrix $C = [p_{i_m'i_k'}] = B - A'$. Observe that $C$ is well-defined and symmetric by the following. 
     \[p_{i_m'i_k'} = b_{i_m'i_k'} - a_{i_m'i_k'} = b_{i_k'i_m'} - a_{i_k'i_m'} = p_{i_k'i_m'} \]
     Since, $ (i_m,i_k\mid i_k,i_m) \rhd A = (i_m',i_k'\mid i_k',i_m') \rhd B$, we have, $b_{i_m'i_k'} - a_{i_mi_k} = \frac{b_{i_m'i_m'}+b_{i_k'i_k'} - a_{i_mi_m} - a_{i_ki_k}}{2}$. Then, the system of equations with variables $p_{i}$, $p_{i}+ p_{j} = p_{ij}$ is consistent with the following solution. 
     \[ p_{i'} = \frac{b_{i'i'}- a'_{i'i'}}{2} =  \frac{b_{i'i'}- a_{ii}}{2}, \text{ for } i= 1,2,\cdots,n\]
Now, define the generalised permutation matrix $P$ as, $P_{ii'} = p_i $, and remaining entries are $\varepsilon$. Then, for this $P$, we have $B = P^tAP$.
 \end{proof}
 \begin{example}
     Let $A = \begin{pmatrix}
         1&3&4&5\\
         3&2&6&7\\
         4&6&-3&8\\
         5&7&8&6
     \end{pmatrix}$, and $B = \begin{pmatrix}
         6&9&6&11\\
         9&6&6&11\\
         6&6&3&8\\
         11&11&8&3
     \end{pmatrix}$. Then, for $P = \begin{pmatrix}
         \varepsilon & \varepsilon &1 &\varepsilon\\
         2&\varepsilon & \varepsilon &\varepsilon\\
         \varepsilon & \varepsilon &\varepsilon &3\\
         \varepsilon &0& \varepsilon &\varepsilon
     \end{pmatrix}$, we have $B = P^tAP$. The corresponding permutation $\sigma = (1,3,4,2)$. Here, under the action of $P$, $1 \mapsto 3$, with a weight $p_1 = 1$, $2 \mapsto 1$, with a weight $p_2 = 2$, $3 \mapsto 4$, with a weight $p_3 = 3$, and $4 \mapsto 2$, with a weight $p_4 = 0$. Then, one can observe that, 
     \begin{align*}  (1,4\mid 2,3) \rhd A = (3+8) - (7+4) = (6+11)- (9+8) =  (3,2\mid 1,4) \rhd B
     \end{align*}
 \end{example}
 \begin{corollary}\label{corr-action}
     Let $A$ and $B$ be $n \times n$ finite symmetric matrices, and $\Psi(A) = \Bar{A}$ and $\Psi(B) = \Bar{B}$. Then, $A \sim B$ if and only if, there exists a $\sigma \in S_n$, such that $(\underline{i},\underline{j})-th$ entry of $\Bar{A}$ is the $(\sigma(\underline{i}),\sigma(\underline{j}))-th$ entry of  $\Bar{B}$.
 \end{corollary}
 \begin{proof}
     The proof is immediate from the previous Theorem \ref{lemma_action}. 
 \end{proof}
With the above Theorem and Corollary, we can add a tropical geometry viewpoint to what we have obtained. In the following subsection, we explore this geometric point of view.
 \subsection{A little geometry}
Here, in this subsection, we will demonstrate that the quotient $\Delta_n(\mathbb{T})/ GL_n(\mathbb{T})$ has a nice geometric significance that is related to the tropical projective space. We will briefly recall the classical projective space and its tropical version for completeness. Readers may go through any manifold theory or algebraic geometry textbooks for a detailed understanding of the topic, such as \cite{tu-manifold}. One can read \cite{tropical-geometry} to understand the geometry in tropical algebra. Recall that, a projective space $\mathbb{P_K}^n$, over a field $\mathbb{K}$, is the set $\mathbb{K}^{n+1}- \{0\}/ \approx$, where $``\approx"$ is the equivalence relation defined as $x \approx y$, if there exists a non-zero $\lambda \in \mathbb{K}$, such that $x = \lambda y$. A point $x$ in $\mathbb{P}^n$ is determined by an equivalence class $[x_1,x_2, \cdots, x_{n+1}]$. In case the field $\mathbb{K}$ is the field of the real or the complex numbers, we have a topology on the projective space given by the quotient topology of the usual topology of $\mathbb{K}^{n+1}$ It also obtains a structure of a manifold with the charts of open sets defined as,
\[U_i = \{[x_1,x_2, \cdots, x_{n+1}] \in \mathbb{P}^n: x_i \neq 0\},\]
with $h_i : U_i \rightarrow \mathbb{K}^{n}$, as $h_i ([x]) = (\frac{x_1}{x_i},\frac{x_2}{x_i}, \cdots,\frac{x_{i-1}}{x_i},\frac{x_{i+1}}{x_i}, \cdots, \frac{x_{n+1}}{x_i})$. Each $\frac{x_k}{x_i}$ is called the projective local coordinates. The collection $\{U_i\}_{i=1}^{n+1}$ gives an open cover to $\mathbb{P}^n$. The local charts satisfy a compatibility condition on the intersections of two charts. This definition of a projective space can be carried forward to the tropical case \cite{tropical-geometry}, with special attention to tropical arithmetic whenever needed.\\
 
 Observe the co-factor $(\underline{i}, \underline{j})$ as a map on the set of matrices to a projective space of matrices.

 Since the set of matrices, especially the set of symmetric matrices, form a tropical vector space, we can define a projective space of such matrices. Namely, the projective space of tropical matrices by identifying $M_n{\mathbb{T}}$ as $\mathbb{T}^{n \times n}$. Consider the subspace symmetric matrices $\Delta_n(\mathbb{T})$ of $ M_n(\mathbb{T})$. Then, the projective space $\mathbb{P}\Delta_n(\mathbb{T}) = \Delta_n(\mathbb{T})-\{\varepsilon\}/\approx $, where $A \approx B$, if $A = \lambda B$, for some $\lambda \in \mathbb{T}^*$, can be viewed as follows. Define the standard charts as
\[U_{ij} = \{[A = (a_{11},a_{12}, \cdots, a_{1j}, \cdots, a_{j1}, \cdots, a_{ij}, \cdots, a_{ji}, \cdots, a_{nn})] \in \mathbb{P}\Delta_n(\mathbb{T}): a_{ij}= a_{ji} \neq \varepsilon\},\]
with $h_{ij}([A]) = (\frac{a_{11}}{a_{ij}},\frac{a_{12}}{a_{ij}}, \cdots,\frac{a_{i(j-1)}}{a_{ij}},\frac{a_{i(j+1)}}{a_{ij}}, \cdots, \frac{a_{nn}}{a_{ij}})$. There are $\binom{n}{2}$ of such open sets, which are isomorphic to the copies of $\Delta_{n-1}(\mathbb{T})$. Then, $\{U_{ij}\}_{i,j}$ gives an open cover for $\mathbb{P}\Delta_n(\mathbb{T})$. Take  $~\mathcal{U} = \cap_{i,j=1}^{n} U_{ij}$. Then, $~\mathcal{U}$ consists of equivalence classes of finite symmetric matrices.\\

 \begin{definition}
     Define $D_{2,n} \subset M_{\binom{n}{2}}(\mathbb{T})$, as the collection of all matrices which satisfy the $Pl\ddot{u}cker$ like relations as defined in \eqref{Plucker-relation.1} and \eqref{Plucker-relation.2}.  
 \end{definition}

As we know if $A= \lambda B$, then their $2 \times 2$ co-factors coincide, thus defining a map, $\Tilde{\pi} : \mathcal{U}  \mapsto D_{2,n}$, as $\Tilde{\pi} = \Psi\mid _{\mathcal{U}}$. Now, in the light of Theorem \ref{lemma_action}, and Corollary \ref{corr-action}, by taking quotient with $S_n$, we have a map $\chi : \mathcal{U} \to  D_{2,n}/ S_n$. Then a fiber on this quotient space $\chi$ is the following. 
 \[{\chi} ^{-1} ([\Bar{A}]) = \{B: B \sim A\}.\] 

 Let us summarise the above in a theorem as follows. 
 \begin{theorem}
     Let $[\Bar{A}]$ be in the image of $[A]$ under the map $\chi : \mathcal{U} \to  D_{2,n}/ S_n$, then the fiber at $[\Bar{A}]$ is the set of symmetric matrices $B$ that are congruent to $A$.
\end{theorem}
     
\subsection{Symmetric matrices} 
\begin{definition}
  Let $A \in M_n(\mathbb{T})$. Then we say $A$ is pseudo-diagonalizable if there exists an invertible matrix $P$, such that $(P^tAP)_{ij} = \begin{cases}
      d_i ~\text{, if } i= j\\
      0 ~\text{, otherwise}
  \end{cases} $, for some $d_i$ $> \varepsilon$.   
\end{definition}

\begin{theorem}\label{diagonalizability}
    Let $A$ be an $n\times n$ finite symmetric matrix, such that, there exists $(d_1,d_2, \cdots, d_n)$, satisfying condition,
    \begin{equation}\label{di_equation}
        d_i +d_j = \begin{vmatrix}
 a_{ii}&a_{ij}\\
 a_{ij}&a_{jj}
 \end{vmatrix}_T =  a_{ii} + a_{jj} -2a_{ij}, ~~ \forall ~i,j.
    \end{equation}
    \[\text{ Then } A \sim D_A := \begin{pmatrix}
        d_1 & 0 &0 &\dots &0\\
        0&d_2&0&\dots & 0\\
       \vdots\\
        0&0&0&\dots & d_n
    \end{pmatrix} \] 
\end{theorem}
\begin{proof}
    Let $P:= \begin{pmatrix}
        \frac{d_1-a_{11}}{2} & \varepsilon &\varepsilon &\dots &\varepsilon\\
        \varepsilon&\frac{d_2-a_{22}}{2} &\varepsilon&\dots & \varepsilon\\
        \vdots\\
        \varepsilon&\varepsilon&\varepsilon&\dots & \frac{d_n-a_{nn}}{2} 
    \end{pmatrix}$. Then, we have $P^t A P = D_A$.\\
    
\end{proof}
\begin{theorem}
    Let $A \in M_n(\mathbb{T})$ be a finite, symmetric matrix. Then $A$ is pseudo-diagonalizable if and only if $ (i,j\mid k,l)\rhd A = 0$, for any quadruple $(i,j,l,k)$ with distinct entries.    
\end{theorem}
\begin{proof}
    Let $A$ be pseudo-diagonalizable. Then we have, $A \sim D_A$, as in above Theorem \ref{diagonalizability}, via $P \in Gl_n(\mathbb{T})$. From Theorem \ref{diagonalizability}, we have seen that the corresponding permutation $\sigma$ for such $P$ is the identity permutation. Then by Theorem \ref{lemma_action}, we have,
    \[ (i,j\mid k,l) \rhd A = (i,j\mid k,l) \rhd D_A. \]
This gives, for distinct $i,j,k,l$, we have $ (i,j\mid k,l) \rhd A = 0$. i.e., any $2 \times 2$ pseudo-minor of $A$, which do not involve any diagonal entry is zero.\\
Now, conversely, assume that for each quadruple $(i,j,k,l)$ with distinct entries, the multi-index operator acts on $A$ gives zero. Then, the coefficient matrix corresponding to the system of equations, $d_i+d_j = a_{ii}+a_{jj}-2a_{ij}$, and the corresponding augmented matrix, both have rank $n$, makes the system consistent. i.e., there exists $d_i$'s which satisfy the system of equations $d_i+d_j = a_{ii}+a_{jj}-2a_{ij}$. Then, by Theorem \ref{diagonalizability}, we have $A$ is pseudo-diagonalizable.   
\end{proof}
\textbf{Remark:} One can find such $d_i$'s by solving the real system, corresponding to the system of equations $d_i +d_j =  a_{ii} + a_{jj} -2a_{ij}$. There will be $\binom{n}{2}$ equations, and $d_1,d_2,\cdots,d_n$ as variables. The coefficient matrix has 2 ones in every row, $n-1$ ones in every column, and the rest are zeros. Note that any finite $2 \times 2$ and $3 \times 3$ symmetric matrices are pseudo-diagonalizable. Moreover, for $3 \times 3$ finite symmetric matrices $A$, the corresponding coefficient matrix of the system of equations is invertible. As a result, the matrix is pseudo-diagonalizable in a unique way. i.e., a unique triplet $(d_1,d_2,d_3)$ exists, such that this triplet forms the diagonal of $D_A$, corresponding to the matrix $A$. 
\begin{corollary}
  Let $A$ be an $n\times n$ finite symmetric, pseudo-diagonalizable matrix.  Then, \[A \sim M_A := \begin{pmatrix}
        0 & \frac{2a_{12}-(a_{11}+a_{22})}{2} &\frac{2a_{13}-(a_{11}+a_{33})}{2}&\dots &\frac{2a_{1n}-(a_{11}+a_{nn})}{2}\\
        \frac{2a_{12}-(a_{11}+a_{22})}{2}&0&\frac{2a_{23}-(a_{22}+a_{33})}{2}&\dots & \frac{2a_{2n}-(a_{22}+a_{nn})}{2}\\
       \vdots\\
        \frac{2a_{1n}-(a_{11}+a_{nn})}{2}&\frac{2a_{2n}-(a_{22}+a_{nn})}{2}&\frac{2a_{3n}-(a_{33}+a_{nn})}{2}&\dots & 0
    \end{pmatrix}\]    
\end{corollary}
\begin{proof}
    Let $D_A$ be the matrix as from the previous Theorem. Then, for $Q := \begin{pmatrix}
        \frac{d_1}{2} & \varepsilon &\varepsilon&\dots &\varepsilon\\
        \varepsilon&\frac{d_2}{2}&\varepsilon&\dots & \varepsilon\\
        \vdots\\
        \varepsilon&\varepsilon&\varepsilon&\dots & \frac{d_n}{2}
    \end{pmatrix}$, we have $Q^tM_AQ = D_A$. But, from the previous Theorem, we have $A \sim D_A$. This gives, $A \sim M_A$. 
\end{proof}
\subsection{Skew-symmetric matrices}
As we mentioned in the introduction, classifying tropical skew-symmetric matrices, as defined in Definition \ref{skew-sym}, is essentially classifying all directed simple graphs. Hence, to hope for a nice structure for its canonical form as in the classical case (for example, the canonical forms described in \cite{darboux1}) is a difficult task. However, we can identify some useful structures. 
\begin{theorem}\label{Skew-can.}
    Let $A$ be an invertible $n \times n$ skew-symmetric matrix, with $n$ being odd. There, there exists permutation matrices $C_n$, and $P$, such that $P^tAP = C_n$. Moreover, the permutations corresponding to $A$ and $C_n$ are conjugate to one another. 
\end{theorem}
\begin{proof}
    Let $A$ be a skew-symmetric permutation matrix. i.e., $A$ is a generalized permutation matrix. Let $\sigma$ be the corresponding permutation to $A$. Note that $\sigma$ can not have any fixed points, as $A$ is skew-symmetric. Let $\tau$ be the corresponding permutation to the desired $C_n$. Now, choose a permutation $\pi_P$ (corresponding to the desired $P$), such that it follows $\pi_P^t \sigma \pi_P = \tau$.\\
    Now, assume that, the finite entries of $A$ are $a_1, a_2, \cdots, a_n$, and $x_1, x_2, \cdots, x_n$ are those of $P$. Then, the finite elements of $P^tAP$ are of the form $a_{\sigma(i)}x_ix_{\tau(i)} $, which has to be $0$. In real terms, $a_{\sigma(i)} + x_i + x_{\tau(i)} = 0$. This gives rise to a real system of equations, with $x_i$'s as variables, and the coefficient matrix has exactly two $1$'s in each row and column, and all other entries are $0$'s. As $n$ is odd, this coefficient matrix is invertible, and hence the system has a unique solution for $(x_1,x_2, \cdots, x_n)^t$. This gives the desired $P$. 
\end{proof}
\begin{corollary}
    Let $A$ be an $n \times n$ skew-symmetric matrix, with $n$ being odd, with the corresponding permutation having a single disjoint cycle of length $n$. Then, there exists a permutation matrix $P$, such that,
    \[P^tAP = C_n = \begin{pmatrix}
        \varepsilon &0&\varepsilon &\cdots & \varepsilon\\
        \varepsilon&\varepsilon&0&\cdots&\varepsilon\\
        &&&\ddots&\\
        \varepsilon&\varepsilon &\varepsilon&\cdots&0\\
        0&\varepsilon&\varepsilon&\cdots&\varepsilon
    \end{pmatrix}.\]
\end{corollary}
\begin{proof}
    The proof is immediate from the above proof, with $\tau = (1,2,3, \cdots,n)$.
\end{proof}
\begin{theorem}
    Let $A$ be a skew-symmetric matrix, with each of its irreducible components having odd dimensions $n_1, n_2, \cdots, n_r$ and invertible. Then, there exists a generalised permutation matrix $P$, such that \[P^tAP = \begin{pmatrix}
        C_{n_1} & & *\\
        & \ddots& \\
        \varepsilon& & C_{n_r}
    \end{pmatrix}.\]
\end{theorem}
\begin{proof}
    Let $Q$ be the permutation matrix which takes $A$ to its Frobenious-Normal form (see \cite{thebook,brualdi}) via $QAQ^t$, and  $A_{n_i}$'s are the irreducible components of $A$. i.e., diagonal blocks in its Frobenious-Normal form. If each of $A_{n_i}$'s are of odd dimension and invertible, then by the above Theorem \ref{Skew-can.}, there exists a generalised permutation matrix $P_i$, and $C_{n_i}$, such that $P_i^tA_{n_i}P = C_{n_i}$. Then, construct the matrix $\Bar{P}$, with $P_i$'s as its diagonal blocks, and the rest are $\varepsilon$. Now, define the matrix $P:= Q^t\Bar{P}$. Clearly, $P$ is a generalized permutation matrix. Then, for this $P$, we have, 
    \[P^tAP = \begin{pmatrix}
        C_{n_1} & & *\\
        & \ddots& \\
        \varepsilon& & C_{n_r}
    \end{pmatrix}.\]
\end{proof}
\section{Applications in solving Generalised eigenvalue problem}
In this section, we describe the solvability of the generalized eigenvalue problems with the use of Canonical forms of symmetric and skew-symmetric matrices described in the previous section. The problem of finding a generalised eigenvalue problem is the following(\cite{Special.gen.eigen,Binding.gen.eigen,Cuninghame.gen.eigen}):\\
For a given pair of square matrices $A, B$, find a scalar $\lambda \in \mathbb{T}$, and a non-trivial vector $x \in \mathbb{T}^n - \{\varepsilon\}$, such that it follows the equation,
\begin{equation}\label{gen.eigen.problem}
    Ax = \lambda Bx
\end{equation}
Collection of all scalars which satisfy the above equation \eqref{gen.eigen.problem} is called generalised eigenvalues, and is denoted by $\Lambda (A,B)$, and collection of such no-trivial vectors are called generalised eigenvectors, and denoted by $V(A,B)$. We say, $(A,B)$ is solvable if $\Lambda (A,B)$ and $V(A,B)$ are non-empty, otherwise called not solvable.
\begin{definition}
    A square matrix $A$ is called nilpotent, if $A^{(k)} = \varepsilon$, for some $k \in \mathbb{N}$.
\end{definition}
\begin{lemma}
    Let $A$ be a finite a finite matrix, and $B$ be nilpotent. Then, $(A,B)$ is not solvable. 
\end{lemma}
\begin{proof}
    Since $A$ is finite, $Ax$ is finite, unless every $x_i$'s are $\varepsilon$. On the other hand, since $B$ is nilpotent, at least one row of $B$ is $\varepsilon$. Hence, $Bx$ contains at least one $\varepsilon$. This makes $(A,B)$ non-solvable. 
\end{proof}
\begin{lemma}
        Let $A$ be a nilpotent matrix and $B$ be any arbitrary matrix. Then either $(A,B)$  or $(A^t,B^t)$ is solvable. 
\end{lemma}
\begin{proof}
    As $A$ is nilpotent, at least one row or column of $A$ must be $\varepsilon$. Without loss of generality, assume the first column of $A$ is $\varepsilon$. For $x = (x_1,x_2, \cdots, x_n)^t$, in the expression $Ax$, $x_1$ does not appear anywhere as first column of $A$ being $\varepsilon$. Then, choosing $x_1 > \varepsilon$, and $x_i = \varepsilon$, for $i=2,3,\cdots, n$, and $\lambda = \varepsilon$, gives $(A,B)$ is solvable. In the case of a row being $\varepsilon$, the above process can be done for $A^t$.  
\end{proof}
Let $A, B$ be symmetric, and $A$ be pseudo-diagonalizable. i.e., There exists a permutation matrix $P$, such that $P^tAP = D_A$, as given in Theorem \ref{diagonalizability}. Then the generalized eigenvalue problem $A x = \lambda Bx$ is equivalent to the following:
\begin{equation*}
    P^tAPP^{-1}x = \lambda P^tBPP^{-1}x\\
    D_A y = \lambda \Bar{B} y
\end{equation*}
where $y = P^{-1}x$ and $\Bar{B} = P^tBP$. Since $B$ is symmetric, $\Bar{B}$ is also symmetric. Hence, it is equivalent to checking the solvability of $(D_A, \Bar{B})$ instead of $(A,B)$. Hence, without loss of generality, to check the solvability of pseudo-diagonalizable symmetric matrices, we can assume one to be in its pseudo-diagonal form. 
\begin{theorem}
    Let $A, B$ be symmetric matrices, and $A$ be in pseudo-diagonal form, with all diagonal entries $d_i$'s being non-positive. Then, a necessary condition for $(A,B)$ to be solvable is the column maximums for each column in $B$ should not come from a single row. 
\end{theorem}
\begin{proof}
    Let $A = \begin{pmatrix}
        d_1 &0& \cdots&0\\
        0&d_2&\cdots&0\\
        \cdots\\
        0&\cdots&0&d_n
    \end{pmatrix}$, $B= \begin{pmatrix}
        b_{11}&b_{12}&\cdots&b_{1n}\\
        b_{12}&b_{22}& \cdots&b{2n}\\
        \cdots\\
        b_{1n}&b_{2n}&\cdots&b_{nn}
    \end{pmatrix}$, and $x = \begin{pmatrix}
        x_1\\
        x_2\\
        \vdots\\
        x_n
    \end{pmatrix}$.\\
    \begin{equation}\label{A^2x}
     \text{ Then, }   A x = \lambda Bx \\
       \implies A^2 x = \lambda AB x 
    \end{equation}
But, $A^2 = \begin{pmatrix}
    0&0&\cdots&0\\
    0&0&\cdots&0\\
    \cdots\\
    0&0&\cdots&0
\end{pmatrix}$, and\\
$\Bar{B} = AB = \begin{pmatrix}
    d_1b_{11}\oplus b_{12} \oplus \cdots \oplus b_{1n}& d_1b_{12}\oplus b_{22} \oplus \cdots \oplus b_{2n}& \cdots& d_1b_{1n}\oplus b_{2n} \oplus \cdots \oplus b_{nn}\\
    b_{11}\oplus d_2 b_{12} \oplus \cdots \oplus b_{1n}& b_{12}\oplus d_2 b_{22} \oplus \cdots \oplus b_{2n}& \cdots& b_{1n}\oplus d_2 b_{2n} \oplus \cdots \oplus b_{nn}\\
    \cdots\\
    b_{11}\oplus b_{12} \oplus \cdots \oplus d_nb_{1n}& b_{12}\oplus b_{22} \oplus \cdots \oplus d_nb_{2n}& \cdots& b_{1n}\oplus b_{2n} \oplus \cdots \oplus d_nb_{nn}\\   
\end{pmatrix}.$\\

By taking $a = (x_1 \oplus x_2 \oplus \cdots \oplus x_n) - \lambda$, we have, $A^2 x = a \begin{pmatrix}
    0\\
    0\\
    \vdots\\
    0
\end{pmatrix} = \hat{a}$. Hence, the equation \eqref{A^2x} becomes,
  \[\Bar{B} x = \hat{a} \]
This system of equations has a solution, if and only if $\Bar{x}$ is a solution of it, where $\Bar{x_j} = (\underset{i \in N}{max}~ \Bar{b_{ij}}-a)^{-1}$. This happens, if and only if the corresponding $N_j$'s makes an exhaustive list of $N$, where, 
\[N_j (\Bar{B}, \hat{a}) = \{ i \in N | \Bar{x_j} = a- \Bar{b_{ij}}\}.\]
Due to the special structure of $\Bar{B}$, and $d_i \leq 0$, we get $N_j= N - \{k\}$, for some $k \in N$.\\
For example, if $b_{11}$ is the dominating element in the first column of $B$, then $N_1 = N -\{1\}$. Similarly, in each column, if the entry from the first row is dominated, we get $N_j =  N -\{1\}$. In general, in $j^{th}$ column, if entry from $k^{th}$ row is dominated, then $N_j = N - \{k\}$. Hence, to get $\cup_{j=1}^n N_j = N$, the maximum element of columns in $B$ should come from at least two distinct rows.  
\end{proof}
The above procedure can be repeated when either of $A, B$ is a skew-symmetric or symmetric matrix with their respective canonical forms.  
\section{Applications to Commuting Matrices}
In this section, we find all symmetric matrices which commute with a given pseudo-diagonalizable finite symmetric matrix. Finding commuting pairs other than their own polynomials brings a lot of applications. Recently, in the field of cryptography, in \cite{crypto}, authors have found usage of such commuting matrices on the tropical congruent transformation of symmetric matrix by circular matrix. Here, we provide an algorithm to find all commuting symmetric matrices of a pseudo-diagonalizable matrix without assuming any particular form to it, such as a circular matrix. If $A$ is a pseudo-diagonalizable matrix, then we are finding complete symmetric solutions to the matrix equation $AX = XA$. Let $P$ be the matrix which takes $A$ to its pseudo-diagonal form $D_A$. From Theorem \ref{diagonalizability}, we have $P$ as a diagonal matrix, and hence symmetric. Assuming $X$ is symmetric, we can transform the problem equivalently to the following.
\begin{align*}
   AX &= XA\\
   PAPP^{-1}XP &= PXP^{-1}PAP \\
   D_A Y &= Y^t D_A \label{D_AY} 
\end{align*}
where, $Y = P^{-1}XP$. The above equation in tropical algebra, in its most general form, for a given matrix $A$, $AX =X^tA$ is an analogous equation of $AX+X^tA=0$ in usual algebra, which is widely studied by many algebraists, notably in \cite{teran,Garcia}. F. De Teran et al., in \cite{teran}, have studied the solution space for the matrix equation $AX+X^tA = 0$, which occurs as a tangent space for the manifold $P^tAP$. It is not very coincidental that this equation occurred here, and the solving technique does not restrict its applications to just finding the commuting symmetric pairs. It will open windows to a large number of fields wherever this equation appears. For the classical linear algebra case, one can go through \cite{teran_2} to understand the historical development of this problem.\\
Let $T$ be the linear transformation that maps $Vec(X) \mapsto Vec(X^t)$. Then equation $D_A Y = Y^t D_A$ is equivalent to,
\[ (I \boxtimes D_A) Vec(Y) = ((D_A \boxtimes I)T) Vec(Y).\]
Let $B = I \boxtimes D_A $ and $ C= (D_A \boxtimes I)T $. Then, the equation is $By = C y$, where $y = Vec(Y)$. We use the algorithm given in \cite{Ax=Bx} to solve this equation. Note that, $T$ is an $ n^2 \times n^2$ permutation matrix. The corresponding permutation is the following.
\begin{align*}
\sigma = (2~(n+1))(3~(2n+1))\cdots (n~(n(n-1)+1))((n+3)~(2n+2)) ((n+4)~(3n+2)) \\
\cdots (2n~(n(n-1)+2))\cdots (n(n-2)+n~ n(n-1)+n-1).
\end{align*}
 The fixed points of the permutation are $kn+k+1$, where $k = 0, \cdots, n-1$, which are corresponding to the diagonal entries. For $D_A$ as given in Theorem \ref{diagonalizability}, $B$ and $C$ have the following block matrix form. 
 \[B= \begin{pmatrix}
     \begin{pmatrix}
         d_1&0&\cdots&0\\
         0&d_2&\cdots&0\\
         \cdots\\
         0&0&\cdots&d_n
     \end{pmatrix}&
     \varepsilon & \cdots &\varepsilon\\
     \varepsilon & \begin{pmatrix}
         d_1&0&\cdots&0\\
         0&d_2&\cdots&0\\
         \cdots\\
         0&0&\cdots&d_n
     \end{pmatrix}& 
     \cdots& \varepsilon\\
     \cdots&\cdots&\cdots&\cdots\\
     \varepsilon&\varepsilon&\cdots& 
     \begin{pmatrix}
         d_1&0&\cdots&0\\
         0&d_2&\cdots&0\\
         \cdots\\
         0&0&\cdots&d_n
     \end{pmatrix}
 \end{pmatrix}\]
\[C = \begin{pmatrix}
    \begin{pmatrix}
        d_1&0&\cdots&0\\
        \varepsilon&\varepsilon&\cdots&\varepsilon\\
        \cdots\\
         \varepsilon&\varepsilon&\cdots&\varepsilon
    \end{pmatrix}&
    \begin{pmatrix}
        \varepsilon&\varepsilon&\cdots&\varepsilon\\
         d_1&0&\cdots&0\\
        \cdots\\
         \varepsilon&\varepsilon&\cdots&\varepsilon
    \end{pmatrix}& \cdots &
    \begin{pmatrix}
       \varepsilon&\varepsilon&\cdots&\varepsilon \\
        \varepsilon&\varepsilon&\cdots&\varepsilon\\
        \cdots\\
         d_1&0&\cdots&0
    \end{pmatrix}\\
    
    \begin{pmatrix}
        0&d_2&\cdots&0\\
        \varepsilon&\varepsilon&\cdots&\varepsilon\\
        \cdots\\
         \varepsilon&\varepsilon&\cdots&\varepsilon
    \end{pmatrix}&
    \begin{pmatrix}
        \varepsilon&\varepsilon&\cdots&\varepsilon\\
         0&d_2&\cdots&0\\
        \cdots\\
         \varepsilon&\varepsilon&\cdots&\varepsilon
    \end{pmatrix}& \cdots &
    \begin{pmatrix}
       \varepsilon&\varepsilon&\cdots&\varepsilon \\
        \varepsilon&\varepsilon&\cdots&\varepsilon\\
        \cdots\\
         0&d_2&\cdots&0
    \end{pmatrix}\\
    \cdots&\cdots&\cdots&\cdots\\
    \begin{pmatrix}
        0&0&\cdots&d_n\\
        \varepsilon&\varepsilon&\cdots&\varepsilon\\
        \cdots\\
         \varepsilon&\varepsilon&\cdots&\varepsilon
    \end{pmatrix}&
    \begin{pmatrix}
        \varepsilon&\varepsilon&\cdots&\varepsilon\\
        0&0&\cdots&d_n\\
        \cdots\\
         \varepsilon&\varepsilon&\cdots&\varepsilon
    \end{pmatrix}& \cdots &
    \begin{pmatrix}
       \varepsilon&\varepsilon&\cdots&\varepsilon \\
        \varepsilon&\varepsilon&\cdots&\varepsilon\\
        \cdots\\
         0&0&\cdots&d_n
    \end{pmatrix} 
\end{pmatrix}\]
Both $B$ and $C$ are in their dominant form with respect to each other already. Hence, their sum matrix $M = B \oplus C$ is the following. 
\[M = \begin{pmatrix}
     \begin{pmatrix}
         d_1&0&\cdots&0\\
         0&d_2&\cdots&0\\
         \cdots\\
         0&0&\cdots&d_n
     \end{pmatrix}&
     \begin{pmatrix}
        \varepsilon&\varepsilon&\cdots&\varepsilon\\
         d_1&0&\cdots&0\\
        \cdots\\
         \varepsilon&\varepsilon&\cdots&\varepsilon
    \end{pmatrix}& \cdots &
    \begin{pmatrix}
       \varepsilon&\varepsilon&\cdots&\varepsilon \\
        \varepsilon&\varepsilon&\cdots&\varepsilon\\
        \cdots\\
         d_1&0&\cdots&0
    \end{pmatrix}\\
     \begin{pmatrix}
        0&d_2&\cdots&0\\
        \varepsilon&\varepsilon&\cdots&\varepsilon\\
        \cdots\\
         \varepsilon&\varepsilon&\cdots&\varepsilon
    \end{pmatrix} & \begin{pmatrix}
         d_1&0&\cdots&0\\
         0&d_2&\cdots&0\\
         \cdots\\
         0&0&\cdots&d_n
     \end{pmatrix}& 
     \cdots& \begin{pmatrix}
       \varepsilon&\varepsilon&\cdots&\varepsilon \\
        \varepsilon&\varepsilon&\cdots&\varepsilon\\
        \cdots\\
         0&d_2&\cdots&0
    \end{pmatrix}\\
     \cdots&\cdots&\cdots&\cdots\\
     \begin{pmatrix}
        0&0&\cdots&d_n\\
        \varepsilon&\varepsilon&\cdots&\varepsilon\\
        \cdots\\
         \varepsilon&\varepsilon&\cdots&\varepsilon
    \end{pmatrix}&
    \begin{pmatrix}
        \varepsilon&\varepsilon&\cdots&\varepsilon\\
        0&0&\cdots&d_n\\
        \cdots\\
         \varepsilon&\varepsilon&\cdots&\varepsilon
    \end{pmatrix}&\cdots& 
     \begin{pmatrix}
         d_1&0&\cdots&0\\
         0&d_2&\cdots&0\\
         \cdots\\
         0&0&\cdots&d_n
     \end{pmatrix}
 \end{pmatrix} \]
Now, the system is equivalent to $By = Cy = My$. To execute the algorithm given in \cite{Ax=Bx}, we need to introduce the following notations.
\begin{definition}[\cite{Ax=Bx}]
    \begin{itemize}
        \item[(i)]  $WB(i) = \{j : b_{ij} > c_{ij}\}, WC(i) = \{j : b_{ij} < c_{ij}\}$.
        \item[(ii)]  $E(i) = \{j : b_{ij} = c_{ij} \neq \varepsilon \}, F(i) = \{j : b_{ij} = c_{ij} =\varepsilon \}.$
        \item[(iii)] $win(i) = (WB(i) \times WC(i)) \cup (E(i) \times E(i)) \subset N \times N$. Each element of $win(i)$ is called a winning pair. 
        
    \end{itemize}
From the above $B,C$ and $M$, we have $win(i)$ as follows.
\begin{align*}
    win(1) &= [n] \times [n]\\
    win(2) &= [n] \times (n +[n])\\
    \cdots\\
    win(n) &= [n] \times (n(n-1)+[n])\\
    win(n+1) &= (n+[n]) \times [n]\\
    win(n+2) &= (n+[n]) \times (n+[n])\\
    \cdots\\
    win (n^2) &= (n(n-1)+[n]) \times (n(n-1)+[n])
    \end{align*}
\end{definition}

   For an $I = (\underline{i}_1,\underline{i}_2) \in win(i)$, by $|I|$ we denote the set $\{\underline{i}_1,\underline{i}_2\}$. For $i, k \in N, \text{ and } i < k $, suppose $I \in win (i)$ and $K \in win(k)$. Then , the value of the $2 \times 2$ tropical minor of $M$, denoted by $M(i, k;\underline{i}, \underline{k})$, is defined as,
   \[\begin{vmatrix}
       m_{i\underline{i}}&m_{i\underline{k}}\\
       m_{k\underline{i}}&m_{k\underline{k}}
   \end{vmatrix}_{trop} = max \{m_{i\underline{i}}+ m_{k\underline{k}}, m_{i\underline{k}}+m_{k\underline{i}}\}, \text{ where } \underline{i} \in |I|, \text{ and } \underline{k} \in |K| \] 
   Readers may note that this tropical determinant is different from the pseudo-determinant defined in Section 3.
\begin{definition}(\cite{Ax=Bx})
     For $i, k \in N, \text{ and } i < k $, let $I \in win (i)$ and $K \in win(k)$. Then we say, $K$ is compatible with $I$ ,if $M(i, k;\underline{i},\underline{k}) = m_{i\underline{i}}+ m_{k\underline{k}} $, for all $\underline{i} \in |I|, \text{ and } \underline{k} \in |K|$. i.e., the tropical determinant corresponding to $I, K$ is attained on main diagonal for all $\underline{i} \in |I|, \text{ and } \underline{k} \in |K|$.
\end{definition}
\begin{definition}(\cite{Ax=Bx})
    Let $\gamma = (I_1,..., I_n)$ be an $n$-tuple with $I_h \in win(h)$, for every $h \in N$. We say that $\gamma$
is a win sequence for $By = Cy = My$, if $I_h$ is compatible with $I_i$, for all $1 \leq i < h \leq m$.
\end{definition}

 Here, for each $k =1,2,\cdots, n^2$, we have,
 \[Win(k) = (\lfloor k/n \rfloor n +N ) \times ((k-1)n+N).\]
 From the structure of $M, B$ and $C$, we can describe all winnable sequences irrespective of the order relation of $ d_i$s. In fact, if we write all winnable sequences starting with a root entry $(i,j)$, then this forms a tree with multiple branches. If any of the branches terminate before having $n$ length walk from the root, that branch can be cut off from the tree of winnable sequences. With this help, we present a modified algorithm given by E. Lorenzo et al. in \cite{Ax=Bx}, which gives all symmetric commuting matrices of a pseudo-diagonalizable matrix.\\ 
 
 \textbf{Algorithm}\\
 \begin{itemize}
     \item[Step 1] Compute the matrix $P$, which (pseudo) diagonalize $A$, via $D_A = P^tAP$.
     \item[Step 2] Find all possible $n-$tuples $(d_1,d_2,\cdots,d_n)$ which fits in the diagonal of $D_A$.
     \item[Step 3] Compute all trees of winnable sequences beginning with the root point $(i,i)$, for each $i \in N$.
     \item[Step 4] For finding winnable sequence tree, with root $(i,j)$ (or $(j,i)$), find all intersections trees with root $(i,i)$ and $(j,j)$.
     \item[Step 5] With respect to each winnable sequence $\gamma$, by using Algorithm in \cite{Ax=Bx}, find $Sol_\gamma = \{y_\gamma:~ By_\gamma=Cy_\gamma\}$.
     \item[Step 6] Find the set $\Bar{\mathcal{Y}} = \{ Y_\gamma : Vec(Y_\gamma) = y_\gamma, y_\gamma \in Sol_\gamma, \text{ for all winnable sequence } \gamma \}.$   
     \item[Step 7] Find $\mathcal{Y} = \{X :~ X = PY_\gamma P^{-1}, Y_\gamma \in \Bar{\mathcal{Y}} \}$.
     \item[Step 8] Find $\mathcal{X} = \{X \in \mathcal{Y}:~ X^t = X\}$. This set, $\mathcal{X}$ gives the set of all symmetric commuting matrices with $A$.
 \end{itemize}
 \begin{example}
   Let $A = \begin{pmatrix}
       8&3\\
       3&8
   \end{pmatrix}$. Then, $|A|_T = 10$. Choose $d_1 =4 $ and $d_2 = 6$. For $P = \begin{pmatrix}
       -2&\varepsilon\\
       \varepsilon&-1
   \end{pmatrix}$, we have $D_A = P^tAP = \begin{pmatrix}
       4&0\\
       0&6
   \end{pmatrix} $. Then, we have the following, 
   \[ B = \begin{pmatrix}
       4&0&\varepsilon&\varepsilon\\
       0&6&\varepsilon&\varepsilon\\
       \varepsilon&\varepsilon&4&0\\
       \varepsilon&\varepsilon&0&6
   \end{pmatrix}, ~ C =\begin{pmatrix}
       4&0&\varepsilon&\varepsilon\\
       \varepsilon&\varepsilon&4&0\\
       0&6&\varepsilon&\varepsilon\\
       \varepsilon&\varepsilon&0&6
   \end{pmatrix}  \text{ and }  M = B \oplus C = \begin{pmatrix}
        4&0&\varepsilon&\varepsilon\\
       0&6&4&0\\
       0&6&4&0\\
       \varepsilon&\varepsilon&0&6
   \end{pmatrix}\]
Then, the trees of winnable sequences are the following. 
 \[ \includegraphics[width=17cm, height=9cm]{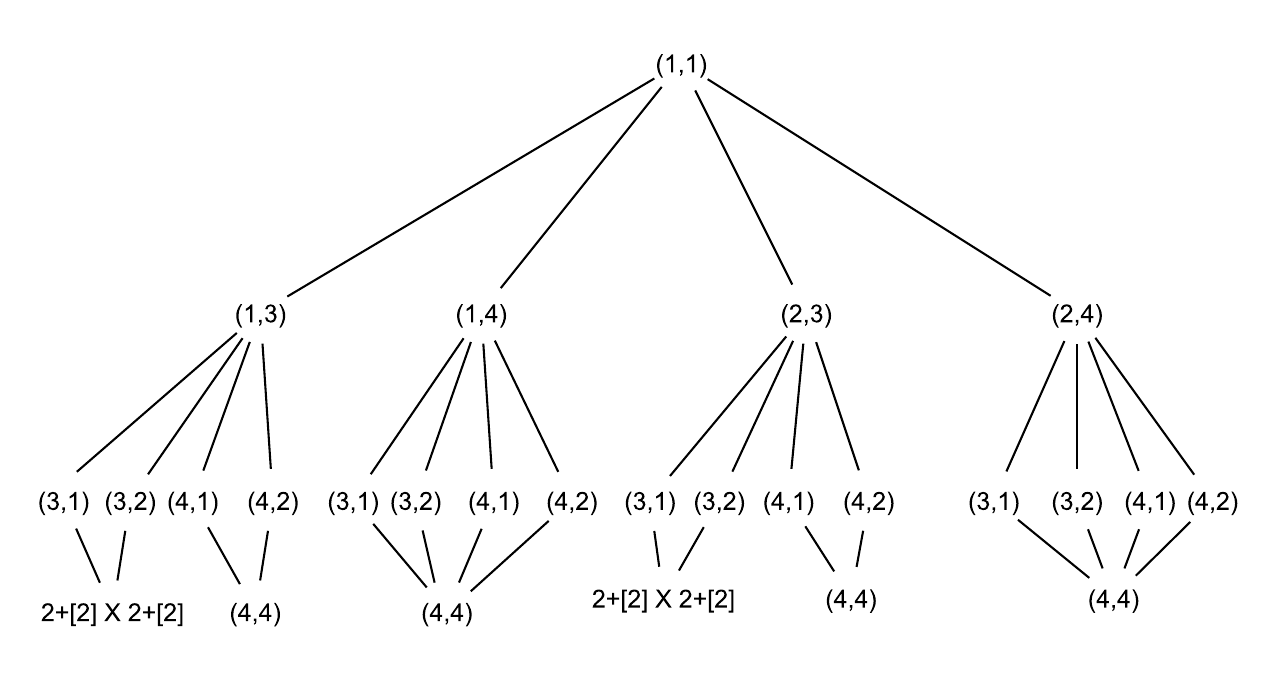}  \]
  \[\begin{tikzpicture}
[
    level 1/.style={sibling distance=55mm},
    level 2/.style={sibling distance=35mm},
]
	\node {$(2,2)$}
		child {node {$(2,3)$}
                child{node {$(3,2)$}
                child{node{$2+[2] \times 2+[2]$}}
                }
                 child{node {$(4,2)$}
                child{node{$2+[2] \times 2+[2]$}}
                }
  }
		child {
		    node {$(2,4)$}
		    child {node {$(3,2)$}
                child {node{$(4,4)$}} 
      }
		    child {node {$(4,2)$}
                child {node{$(4,4)$}} 
      }
		};
\end{tikzpicture}\]
In the tree, each level corresponds to each index $i$, top to bottom in increasing order, and  $[n]$ represents $\{1,2,\cdots,n\}$. Note that, for any $2 \times 2$ $D_A$, with $d_i \geq 0$, the above two trees give an exhaustive list of winnable sequences. Now, let us choose a winnable sequence and construct a solution. Take
\[\gamma =~ ((1,2),(2,3),(3,2),(4,4))\]
Let $Y = \begin{pmatrix}
    x_1&x_3\\
    x_2&x_4
\end{pmatrix}$, so that $y = Vec(Y) = \begin{pmatrix}
    x_1\\
    x_2\\
    x_3\\
    x_4
\end{pmatrix}$.
Then, the corresponding set of equations and inequalities are the following. 
\begin{align*}
   \begin{cases}
       x_2 -x_1 &= 4\\
    x_3 - x_2 &= 2
   \end{cases}\\
   \begin{cases}
        x_1 &\leq 6+x_2\\
        x_4 &\leq 6+ x_2\\
        x_1 &\leq 4+x_3\\
        x_4 &\leq 4+x_3\\
        x_3 &\leq 6+x_4
   \end{cases}   
\end{align*}
On solving this, we get, $\mathcal{\Bar{Y}_\gamma} = \Big\{ \begin{pmatrix}
    x_1&6+x_1\\
    4+x_1&x_4
    \end{pmatrix}: x_4-x_1 \leq 10, ~ x_1 \leq x_4 \Big\}$. Now, applying the action of the matrix $P$, we get,
    \[\mathcal{Y_\gamma} = P\mathcal{\Bar{Y}_\gamma}P^{-1} =  \Big\{ \begin{pmatrix}
    x_1&5+x_1\\
    5+x_1&x_4
    \end{pmatrix}: x_4-x_1 \leq 10, ~ x_1 \leq x_4 \Big\}.\]
i.e., the parameters $x_1,x_4$ varies in the following region.\\
\[\begin{tikzpicture}
    \begin{axis}[
        axis lines = middle,
        xlabel = {$x_1$},
        ylabel = {$x_4$},
        xmin=-10, xmax=20,
        ymin=-10, ymax=30,
        xtick={-10,-5,0,5,10,15,20},
        ytick={-10,-5,0,5,10,15,20,25,30},
        legend pos=north west
    ]
    
    \addplot[
        name path=upper, 
        domain=-10:20, 
        samples=100, 
        color=black,
    ]
    {x + 10};
    
    \addplot[
        name path=lower, 
        domain=-10:20, 
        samples=100, 
        color=black,
    ]
    {x};
    
    \addplot [
        thick,
         fill opacity=0.7, 
        color=gray!30,
    ] fill between[
        of=upper and lower,
    ];
    
    \end{axis}
\end{tikzpicture}\]
Since this set is completely symmetric, we can take all of them into $\mathcal{X}_\gamma$. Then, for any $X \in \mathcal{X}_\gamma$, we have,
\[ \begin{pmatrix}
    8&3\\
    3&8
\end{pmatrix} \begin{pmatrix}
    x_1&5+x_1\\
    5+x_1&x_4
\end{pmatrix} = \begin{pmatrix}
     x_1 +8 & x_1+13\\
     x_1+13&x_4+8
\end{pmatrix} =  \begin{pmatrix}
    x_1&5+x_1\\
    5+x_1&x_4
\end{pmatrix} \begin{pmatrix}
    8&3\\
    3&8
\end{pmatrix}. \]
Similarly, for each winnable sequence $\gamma$, we can generate a set of symmetric commuting matrices. Also, note that every winnable sequence does not need to generate symmetric matrices, hence a solution to our problem. But this method guarantees to find every solution to the problem of finding symmetric commuting matrices with $A$. 
    \end{example}
\bibliographystyle{abbrv}
\bibliography{monoid}
 \end{document}